\documentclass[10pt]{article}
\usepackage{amsfonts}
\usepackage{amsmath}
\usepackage{mathrsfs,multirow}
\usepackage{mathrsfs,amscd,amssymb,amsthm,amsmath,bm,graphicx,psfrag,subfigure}

\makeatletter

\renewcommand{\@seccntformat}[1]{{\csname the#1\endcsname}{\normalsize .}\hspace{.5em}}
\makeatother

\def \[{\begin{equation}}
\def \]{\end{equation}}

\newtheorem{thm}{Theorem}[section]

\newtheorem{lem}[thm]{Lemma}

\begin{document}
\setlength{\baselineskip}{15pt}
\begin{center}{\Large\bf The extremal problems on the inertia of weighted bicyclic graphs\footnote{Financially supported by the National Natural Science Foundation of China (Grant Nos. 11271149, 11371062), the Program for New Century Excellent Talents in University (Grant No. NCET-13-0817) and the Special Fund for Basic Scientific Research of Central Colleges (Grant No. CCNU13F020)).}}

\vspace{4mm}
{\large Shibing Deng, Shuchao Li\footnote{Corresponding author. \\
\hspace*{5mm}{\it Email addresses}: 750861119@qq.com (S.B. Deng), lscmath@mail.ccnu.edu.cn (S.C.
Li), 928046810@qq.com (F.F. Song)}, Feifei Song}\\[5pt]

Faculty of Mathematics and Statistics,  Central China Normal
University, Wuhan 430079, P.R. China
\end{center}
\noindent {\bf Abstract}:
Let $G_w$ be a weighted graph. The number of the positive, negative and zero eigenvalues in the spectrum of $G_w$ are called positive inertia index, negative inertia index and nullity of $G_w$, and denoted by $i_{+}(G_w)$, $i_{-}(G_w)$, $i_{0}(G_w)$, respectively. In this paper, sharp lower bound on the positive (resp. negative) inertia index of weighted bicyclic graphs of order $n$ with pendant vertices is obtained. Moreover, all the weighted bicyclic graphs of order $n$ with at most two positive, two negative and at least $n-4$ zero eigenvalues are identified, respectively.

\vspace{2mm} \noindent{\it Keywords}: Weighted bicyclic graphs; Adjacency matrix; Inertia

\vspace{2mm}

\noindent{AMS subject classification:} 05C50;\ 15A18

 {\setcounter{section}{0}
\section{\normalsize Introduction}\setcounter{equation}{0}
In this paper, we only consider simple weighted graphs on positive weight set. Let $G_w$ be a weighted graph with vertex set $\{v_{1}, v_{2}, \ldots , v_{n}\}$, edge set $E(G)\neq \emptyset$ and $W(G_w)=\{w_{j}> 0, j=1, 2, \ldots , |E(G)|\}$. The function $w: E(G) \rightarrow W(G_w)$ is called a weight function of $G_w$. It is obvious that each weighted graph corresponds to a weight function. The \textit{adjacency matrix} of $G_w$ on $n$ vertices is defined as the matrix $A(G_w)=(a_{ij})$ such that $a_{ij}=w(v_iv_j)$ if $v_iv_j\in E(G)$ and 0 otherwise. The eigenvalues $\lambda_1,\lambda_2,\ldots,\lambda_n$ of $A(G_w)$ are said to be the eigenvalues of the weighted graph $G_w$. The \textit{inertia} of $G_w$ is defined to be the triple $In(G_w)=(i_+(G_w),i_-(G_w),i_0(G_w))$, where $i_+(G_w),i_-(G_w)$ and $i_0(G_w)$ are the numbers of the positive, negative and zero eigenvalues of $A(G_w)$ including multiplicities, respectively.  $i_+(G_w)$ and $i_-(G_w)$ are called the \textit{positive, negative index of inertia} (for short, \textit{positive, negative index}) of $G_w$, respectively. The number $i_0(G_w)$ is called the \textit{nullity} of $G_w$. Obviously, $i_+(G_w)+i_-(G_w)+i_0(G_w)=n$.

An \textit{induced subgraph} of $G_w$ is an induced subgraph of $G$ having the same weights with those of $G_w$. For an induced weighted subgraph $H_w$ of the weighted graph $G_w$, let $G_w-H_w$ be the subgraph obtained from $G_w$ by deleting all vertices of $H_w$ and all incident edges. We define that the union of $G_{w}^{1}$ and $G_{w}^{2}$, denoted by $G_{w}^{1}\bigcup G_{w}^{2}$, is the graph with vertex-set $V(G_{w}^{1})\bigcup V(G_{w}^{2})$, edge-set $E(G_{w}^{1})\bigcup E(G_{w}^{2})$ and the weight of each edge is not changed. A \textit{bicyclic graph} is a simple connected graph in which the number of edges equals the number of vertices plus 1. A weighted path and a weighted cycle of order $n$ are denoted by $P^{n}_w$, $C^{n}_w$, respectively. An isolated vertex is sometime denoted by $K_{1}$.

The study of eigenvalues of graph has been received a lot of attention due to its applications in chemitry (see \cite{D-2,D-3,D-5,D-8} for details). As we know, if $G$ is a bipartite graph, then $i_+(G)=i_-(G)=\alpha(G)=\frac{n-i_0(G)}{2}$, where $\alpha(G)$ is the matching number of $G$, otherwise, $i_+(G),i_-(G)$ and $i_0(G)$ do dot have this relationship. Gregory et al. \cite{D-6} studied the subadditivity of the positive, negative indices of inertia and developed certain properties of Hermitian rank which were used to characterize the biclique decomposition number. Gregory et al. \cite{D-7} investigated the inertia of a partial join of two graphs and established a few relations between inertia and biclique decompositions of partial joins of graphs. Daugherty \cite{D-9} characterized the inertia of unicyclic graphs in terms of matching number and obtained a linear-time algorithm for computing it. Yu et al. \cite{D-4} investigated the minimal positive index of inertia among all unweighted bicyclic graphs of order $n$ with pendants, and characterized the bicyclic graphs with positive index 1 or 2. Very recently, it is interesting to see that Marina et al. \cite{D-17} studied the inertia set of a signed graph in algebraic approach.

The nullity of unweighted graphs has been studied well in the literature. Tan and Liu \cite{D-11} gave the nullity set of unicyclic graphs and characterized the unicyclic graphs with maximum nullity. In addition, Nath and Sarma \cite{D-12} presented another version of characterization of an acyclic or unicyclic graph to be singular. One of the present authors \cite{D-10} investigated the nullity of graphs with pendant vertices. Fan and Qian \cite{D-15} characterized the bipartite graphs with the second largest nullity and the regular bipartite graphs with the third largest nullity. Fan and Wang \cite{D-14} characterized the unicyclic signed graphs of order $n$ with nullity $n-2,n-3,n-4,n-5$, respectively.

Our paper is motivated directly by \cite{D-1,D-13,D-10,D-18}. On the one hand, Fan et al. \cite{D-1} studied the nullity of signed bicyclic graph (which is, in fact, the bicyclic graph with edge weight 1 or $-1$); Li \cite{D-10} and Hu \cite{D-13} studied the nullity of unweighted bicyclic graph. On the other hand, Yu et al. \cite{D-4} characterized all $n$-vertex unweighted bicyclic graphs with positive index 1 or 2. It is natural and interesting for us to consider the extremal problems on the inertia of weighted bicyclic graphs, which may generalize corresponding results of \cite{D-1,D-13,D-10,D-4}.

This paper is organized as follows: in Section 2, some preliminaries are introduced. In Section 3, we present the lower bound for the positive, negative index of $n$-vertex weighted bicyclic graphs with pendants. In Section 4, we characterize all $n$-vertex weighted bicyclic graphs without pendant twins having one or two positive (resp. negative) eigenvalues. In Section 5, we characterize all $n$-vertex weighted bicyclic graphs without pendant twins of rank $2,3,4$.

\section{\normalsize Preliminaries}\setcounter{equation}{0}

In this section, we list some lemmas which will be used to prove our main results. Suppose $M$, $N$ are two Hermitian matrices of order $n$, if there exists an invertible matrix $Q$ of order $n$ such that $QMQ^*=N$, $Q^*$ denotes the conjugate transpose of $Q$, then we say that $M$ is \textit{congruent to} $N$, denoted by $M\cong N$.
\begin{lem}[\cite{D-16}]\label{lem2.1}
Let $M,N$ be two Hermitian matrices of order $n$ such that $M\cong N$, then $i_+(M)=i_+(N), i_-(M)=i_-(N)$ and $i_0(M)=i_0(N)$.
\end{lem}

It is easy to obtain the following result.
\begin{lem}\label{lem2.2}
Let $G_w=G_{w}^{1}\bigcup G_{w}^{2}\bigcup \ldots \bigcup G^t_{w}$ be a weighted graph, where $G^i_{w}$ $(i=1,2,\ldots, t)$ are connected components of $G_w$. Then $i_{+}(G_w)=\sum^t_{i=1}i_{+}(G^i_{w}), i_{-}(G_w)=\sum^t_{i=1}i_{-}(G^i_{w})$ and $i_{0}(G_w)=\sum^t_{i=1}i_{0}(G^i_{w})$.
\end{lem}

Let $M$ be a Hermitian matrix. We denoted three types of elementary congruence matrix operations (ECMOs) on $M$ as follows:
\begin{enumerate}
                                                           \item  interchanging $i$th and $j$th rows of $M$, while interchanging $i$th and $j$th columns of $M$;                 \item multiplying $i$th row of $M$ by a non-zero number $k$, while multiplying $i$th column of $M$ by $k$;
                                                           \item adding $i$th row of $M$ multiplied by a non-zero number $k$ to $j$th row, while adding $i$th column of $M$ multiplied by $k$ to $j$th column.
\end{enumerate}
By Lemma 2.1, the ECMOs do not change the inertia of a Hermitian matrix.
\begin{lem}[\cite{D-4}]
Let $M$ be an $n\times n$ Hermitian matrix and $N$ be the Hermitian matrix obtained by bordering $M$ as followings:
\begin{eqnarray*}
  N &=& \left(
          \begin{array}{cc}
            M & y \\
            y^* & a \\
          \end{array}
        \right),
\end{eqnarray*}
where $y$ is a column vector, $y^*$ denotes the conjugate transpose of $y$ and $a$ is a real number. Then $i_+(N)-1 \leq i_+(M)\leq i_+(N), i_-(M)-1 \leq i_-(M)\leq i_{-}(N)$.
\end{lem}

By Lemma 2.3  we can get the following result immediately:
\begin{lem}[\cite{D-189}]
Let $H_w$ be an induced subgraph of G. Then $i_+(H_w)\leq i_+(G_w)$ and $i_-(H_w)\leq i_-(G_w)$.
\end{lem}
\begin{lem}[\cite{D-189}]
Let $C_{w}^{n}$ be a weighted cycle of order $n$. Then
$$
i_{+} (C_{w}^{n})= \left\{
                     \begin{array}{ll}
                       \frac{n+1}{2}, & \hbox{if $n\equiv1 \pmod{4};$} \\
                       \frac{n}{2}, & \hbox{if $n\equiv2 \pmod{4};$} \\
                       \frac{n-1}{2}, & \hbox{if $n\equiv3 \pmod{4}.$}
                     \end{array}
                   \right.
i_{-} (C_{w}^{n})= \left\{
                     \begin{array}{ll}
                       \frac{n-1}{2}, & \hbox{if $n\equiv1 \pmod{4};$} \\
                       \frac{n}{2}, & \hbox{if $n\equiv2 \pmod{4};$} \\
                       \frac{n+1}{2}, & \hbox{if $n\equiv3 \pmod{4}.$}
                     \end{array}
                   \right.
$$
Furthermore, if $n\equiv 0 \pmod{4},$ let $C_{w}^{n}=v_{1}v_{2}\ldots v_{n}v_{1}$ be a weighted cycle of order $n$, $w(v_{i}v_{i+1})=a_{i}$ $(1\leq i\leq n)$ and let $v_{n+1}=v_{1}$. Then
$$
i_{+} (C_{w}^{n})=i_{-} (C_{w}^{n})=\left\{
  \begin{array}{ll}
    \frac{n}{2}-1, & \hbox{if $\prod_{i=1}^{\frac{n}{2}}a_{2i-1}= \prod_{i=1}^{\frac{n}{2}}a_{2i};$} \\[5pt]
    \frac{n}{2}, & \hbox{otherwise.}
  \end{array}
\right.
$$
\end{lem}
\begin{lem}[\cite{D-189}]
Let $G_w$ be a graph containing a pendant vertex $v$ with its unique neighbor $u$. Then $i_+(G_w) =i_+(G_w-u-v)+1$ and $i_-(G_w) =i_-(G_w-u-v)+1.$
\end{lem}

The following result is an immediate consequence of Lemma 2.6.
\begin{lem} Let $P^n_w$ be a weighted path of order $n$. Then
$$
i_+(P^n_w) =i_-(P^n_w)=\left\{
                \begin{array}{ll}
                  \frac{n-1}{2}, & \hbox{if $n$ are odd;} \\[3pt]
                  \frac{n}{2}, & \hbox{if $n$ are even;}
                \end{array}
              \right.
$$
\end{lem}
Let $u,v$ be two pendant vertices of a weighted graph $G_w$, $u,v$ are called a \textit{pendant twin} if they have the same neighborhood in $G_w$.
The following result is an immediate consequence of Lemma 2.6 since $i_+(K_1) = i_-(K_1) = 0$.
\begin{lem}
If $u,v$ is a pendant twin in a weighted graph $G_w$, then $i_+(G_w) = i_+(G_w-v) = i_+(G_w-u)$ and $i_-(G_w) = i_-(G_w-v) = i_-(G_w-u).$
\end{lem}

Let $S^k_w$ be a weighted star of order $k$ with center $v$ and non-central vertices $v_1,\ldots, v_{k-1}$. We can get the following two lemmas by Lemmas 2.4 and 2.6.
\begin{lem}
Let $G_{w}^{0}$ be a weighted graph of order $n-k$ such that $u\in V(G_{w}^{0})$. Let $G_{w}^{1}$ be the graph obtained from $G_{w}^{0}$ and $S_w^{k}$ by inserting an edge between $u$ and the center $v$ of $S^k_w$. Let $G_{w}^{2}=G_{w}^{1}-\{vv_1, vv_2,\ldots,vv_{k-1}\} +\{uv_1,uv_2,\ldots, uv_{k-1}\}$ where $w(uv_i)=w(vv_i)$. Then $i_+(G_{w}^{1})\geq i_+(G_{w}^{2})$ and $i_-(G_{w}^{1})\geq i_-(G_{w}^{2})$.
\end{lem}
\begin{lem}
Let $G_{w}^{0}$ be a weighted graph of order $n-l-t$ and $u_1,u_2\in V(G_0)$. Assume that $G_{w}^{1}$ is the graph obtained from $G_{w}^{0}$, $S^{l+1}_w$ and $S^{t+1}_w$ by identifying $u_1$ with the center of $S^{l+1}_w$, $u_2$ with the center of $S^{t+1}_w$,respectively. Let $G^2_w$ be the graph obtained from $G_{w}^{0}$, $S^{l+t+1}_w$ by identifying $u_1$ with the center of $S_{l+t+1}$. Then $i_+(G_{w}^{1})\geq i_+(G_{w}^{2})$ and $i_-(G_{w}^{1})\geq i_-(G_{w}^{2})$.
\end{lem}
\begin{lem}
Let $G_{w}^{1}$ and $G_{w}^{2}$ be two weighted graphs with $u\in V(G_1)$ and $v\in V(G_2)$. Let $P^l_w(l\geq 3)$ be a weighted path with two end-vertices $v_1$, $v_l$. Let  $S^{l}_w$ be a weighted star of order $l$ and have the same weight set with $P^l_w(l\geq 3)$. Let $G'_w$ be the graph obtained from $G^1_w \bigcup G^2_w \bigcup P^l_w$ by identifying $u$ with $v_1$ and $v$ with $v_l$, respectively. Let $G''_w$ be the graph obtained from $G^1_w \bigcup G^2_w$ by identifying $u$, $v$ with the center of $S^{l}_w$. Then $i_+(G'_w)\geq i_+(G''_w)$ and $i_-(G'_w)\geq i_-(G''_w)$.
\end{lem}
\begin{proof} In view of Lemma 2.6, we have
$$
i_+(G''_w)=1+i_+(G_{w}^{1}-u)+i_+(G_{w}^{2}-v).
$$
Note that $(G_{w}^{1}-u)\bigcup (G_{w}^{2}-v)\bigcup P^{l-1}_w$ is an induced subgraph of $G'_w$. In light of Lemma 2.3, it follows that
$$
i_+(G_{w}^{1}-v)+i_+(G_{w}^{2}-u)+i_+(P^{l-1}_w)\leq i_+(G'_w).
$$
By Lemma 2.7, $i_+(P^{l-1}_w)\geq 1\, (l\geq 3)$, therefore $i_+(G'_w)\geq i_+(G''_w)$. Similarly, $i_-(G'_w)\geq i_-(G''_w)$, as desired.
\end{proof}

\section{\normalsize The minimal positive (negative) index of inertia of weighted bicyclic graphs}\setcounter{equation}{0}
 Let $G$ be a bicyclic graph. The \textit{base} of $G$, denoted by $\chi(G)$, is the unique bicyclic subgraph of $G$ containing no pendant vertices. Thus $G$ can be obtained from $\chi(G)$ by attaching trees to some vertices of $\chi(G)$. Let $C^p (p\geq 3)$ and $C^q (q\geq 3)$ be two vertex-disjoint cycles of length $p,q$ and $P_l = v_1v_2\ldots v_l\, (l\geq 1)$ be a path of length $l-1$. Assume that $v\in V(C^p)$ and $u\in V(C^q)$. Let $\infty(p,l,q)$ be the graph obtained from $C^p,C^q$ and $P_l$ by identifying $v$ with $v_1$, $u$ with $v_l$. Let $P_{p+2}, P_{l+2}, P_{q+2}$ be three paths with $\min\{p,l,q\}\geq 0$ and at most one of $p,l,q$ is 0. Let $\theta(p,l,q)$ be the graph obtained from $P_{p+2}, P_{l+2}$ and $P_{q+2}$ by identifying the three initial vertices and terminal vertices. The weighted graphs $\infty(p,l,q)_w$ and $\theta(p,l,q)_w$ are depicted in Fig. 1, where the number on each edge denotes its weight. In what follows in our context, we always assume that the weight for each edge of $\infty(p,l,q)_w$ (resp. $\theta(p,l,q)_w$) are as shown in Fig. 1.
\begin{figure}[h!]
\begin{center}
  \psfrag{2}{$C^p$}\psfrag{1}{$v_1$}
  \psfrag{3}{$C^q$}\psfrag{S}{$b_i$}
  \psfrag{5}{$P_{l+2}$}\psfrag{8}{$v_l$}
  \psfrag{4}{$P_{p+2}$}
  \psfrag{6}{$P_{q+2}$}
  \psfrag{T}{$a_1$}\psfrag{V}{$a_i$}
  \psfrag{U}{$a_p$}\psfrag{M}{$c_1$}
  \psfrag{N}{$c_2$}\psfrag{y}{$c_{l-1}$}
  \psfrag{Q}{$b_1$}\psfrag{R}{$b_{q}$}
  \psfrag{A}{$a_1$}\psfrag{B}{$a_2$}\psfrag{C}{$a_p$}
  \psfrag{D}{$a_{p+1}$}\psfrag{E}{$b_1$}
  \psfrag{F}{$b_2$}\psfrag{P}{$b_{l+1}$}
  \psfrag{I}{$c_1$}\psfrag{J}{$c_2$}\psfrag{K}{$c_q$}
  \psfrag{7}{$P_l$}\psfrag{L}{$c_{q+1}$}\psfrag{f}{$u$}\psfrag{k}{$v$}
  \psfrag{a}{$\infty(p,l,q)_w$}\psfrag{b}{$\theta(p,l,q)_w$}
   \includegraphics[width=120mm]{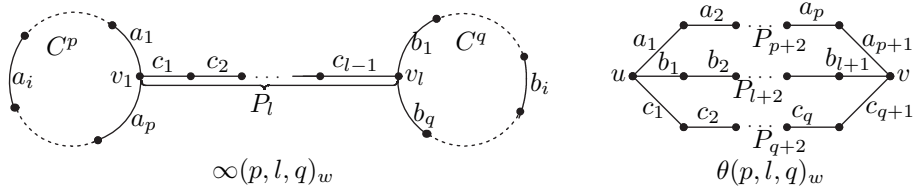}\\
  \caption{Weighted graphs $\infty(p,l,q)_w$ and $\theta(p,l,q)_w.$}
\end{center}
\end{figure}

As we know, the connected bicyclic graphs can be partitioned into two classes: one class of bicyclic graphs contain $\infty(p,l,q)$ as its basis and the other class of bicyclic graphs contain $\theta(p,l,q)$ as its basis. We call bicyclic graph $G$ an $\infty$-graph if $G$ contains some $\infty(p,l,q)$ as its basis and a $\theta$-graph if $G$ contains some $\theta(p,l,q)$ as its basis. We denote by $\mathscr{B}$ (resp. $\mathscr{B}_p$) the set of all weighted bicyclic graphs (resp. weighted bicyclic graphs with pendants) of order $n$. 
Let $\chi(G_w)$ be the base of $G_w$, by Lemma 2.6, there is no correlation between the inertia index of $G_w$ and the weighted set of $G_w-\chi(G_w)$. Hence, in order to determine $In(G_w)$, it suffices to consider the weight of $\chi(G_w)$ in what follows.
\begin{thm} Let $G_w\in {\mathscr{B}_p}$ and contain $\infty(p,l,q)$ as its base. Then
$$
i_+(G_w) \geq \left\{
                \begin{array}{ll}
                  \frac{p+q}{2}, & \hbox{if $p,q$ are odd;} \\[3pt]
                  \frac{p+q}{2}-1, & \hbox{if $p,q$ are even;} \\[3pt]
                  \frac{p+q-1}{2}, & \hbox{otherwise.}
                \end{array}
              \right.
$$
This bound is sharp.
\end{thm}
\begin{proof}
For a weighted $\infty$-graph, let $u$ be the common vertex of $C^p$ and $C^q$ in $\infty(p,1,q)$. Let $G^*$ be the bicyclic graph obtained by attaching $n-p-q+1\,(n\geq p+q)$ pendants to $u$ (see Fig. 2) and let $G^{*}_w$ denote the weighted graph with $G^{*}$ as its underlying graph.
\begin{figure}[h!]
\begin{center}
  \psfrag{1}{$C^p$}\psfrag{6}{$u$}
  \psfrag{2}{$C^q$}\psfrag{7}{$v$}
  \psfrag{8}{$P_{p+2}$}\psfrag{9}{$P_{l+2}$}
  \psfrag{5}{$v_2$}\psfrag{0}{$P_{q+2}$}
  \psfrag{a}{$\infty(p,2,q)$}\psfrag{b}{$G^*$}
  \psfrag{c}{$G^{**}$}
   \includegraphics[width=80mm]{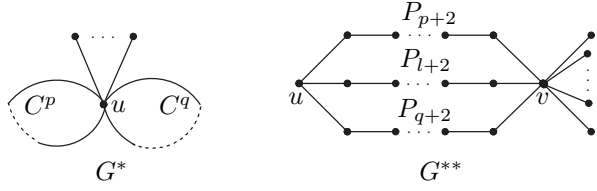}\\
  \caption{Graphs $G^*$ and $G^{**}.$}
\end{center}
\end{figure}

Let $\mathscr{S}_1(n)$ be the set of all $n$-vertex weighted bicyclic graphs whose underlying graph is obtained from $\infty(p,1,q)$ by attaching $n-p-q+1\,(n\geq p+q)$ pendants to a vertex, different from $u$ of $C^p$ or $C^q$.

Let $\mathscr{S}_2(n)$ be the set of all $n$-vertex weighted bicyclic graphs whose underlying graph is obtained from $\infty(p,2,q)$ by attaching $n-p-q\,(n\geq p+q+1)$ pendants to the vertex $u$ of $C^p$ or $C^q$.

In view of Lemma 2.6, we have
\[
i_+(G^*_w) = 1+i_+(P^{p-1}_w)+i_+(P^{q-1}_w).
\]

First we are to show that $i+(G^*_w)\leq i_+(G_w)$ for any $G_w\in \mathscr{S}_1(n)\bigcup \mathscr{S}_2(n).$
In fact, if $G_w\in \mathscr{S}_1(n)$, without loss of generality, we suppose all the pendant vertices are attached at $C^p$. Then by Lemma 2.6, we have
\begin{eqnarray}
i_+(G_w) &=&1+ \left\{
                 \begin{array}{ll}
                   \frac{p-1}{2}+i_+(P^{q-1}_w), & \hbox{if $p$ is odd;} \\[3pt]
                   \frac{p-2}{2}+i_+(C^{q}_w)\ or,\ i_+(P^{q-1}_w), & \hbox{if $p$ is even.}
                 \end{array}
               \right.\notag
   \\
   &=& 1+i_+(P^{p-1}_w)+\left\{
                        \begin{array}{ll}
                          i_+(P^{q-1}_w), & \hbox{if $p$ is odd;} \\[3pt]
                        i_+(C^{q}_w)\ or,\ i_+(P^{q-1}_w), & \hbox{if $p$ is even.}
                        \end{array}
                      \right.
\end{eqnarray}
By Lemma 2.4, $i_+(C^{q}_w)\geq i_+(P^{q-1}_w)$. Hence, in view of (3.1) and (3.2) we have $i_+(G^*_w)\leq i_+(G_w)$.

If $G_w\in \mathscr{S}_2(n)$, without loss of generality, we suppose all the pendant vertices are attached at $C^p$. Then by Lemma 2.6 we have
\begin{eqnarray}
  i_+(G_w) &=&1+ \left\{
                 \begin{array}{ll}
                   \frac{p-1}{2}+i_+(C^{q}_w), & \hbox{if $p$ is odd;} \\[3pt]
                  \frac{p-2}{2}+i_+(C^{q}_w) \ or\ i_+(G'_w), & \hbox{if $p$ is even.}
                 \end{array}
               \right.\notag
   \\
   &=& 1+i_+(P^{p-1}_w)+\left\{
                           \begin{array}{ll}
                             i_+(C^{q}_w), & \hbox{if $p$ is odd;} \\[3pt]
                           i_+(C^{q}_w)\ or\ i_+(G'_w), & \hbox{if $p$ is even,}
                           \end{array}
                         \right.
\end{eqnarray}
where $G'_w$ is a graph obtained by attaching a pendant vertex to a vertex of $C^q$. Note that $i_+(C^{q}_w)\geq i_+(P^{q-1}_w)$ and $i_+(G'_w)\geq i_+(P^{q-1}_w)$ from Lemma 2.4. Hence, in view of (3.1) and (3.3) we have $i_+(G^*_w)\leq i_+(G_w)$.

From Lemmas 2.9, 2.10 and 2.11, $G^*_w$ attains the minimal positive index among all $n$-vertex weighted bicyclic graphs with pendant vertices containing two edge disjoint weighted cycles $C^p_w$ and $C^{q}_w$.
\end{proof}

Similarly, we can have the following theorem:
\begin{thm}
Let $G_w\in {\mathscr{B}_p}$ and contain $\infty(p,l,q)$ as its base. Then
\begin{eqnarray*}
  i_-(G_w) &\geq & \left\{
                     \begin{array}{ll}
                      \frac{p+q}{2}, & \hbox{if $p,q$ are odd;} \\[3pt]
                     \frac{p+q}{2}-1, & \hbox{if $p,q$ are even;} \\[3pt]
                     \frac{p+q-1}{2}, & \hbox{otherwise.}
                     \end{array}
                   \right.
\end{eqnarray*}
This bound is sharp.
\end{thm}

By Theorems 3.1 and 3.2, it follows that
\begin{thm}
Let $G_{w}$ be a weighted $\infty$-graph of order $n$ with pendant vertices, then $i_+(G_w)\geq 3$, $i_-(G_w)\geq 3$ and $i_0(G_w)\leq n-6$.
\end{thm}
\begin{thm}
Let $G_w\in {\mathscr{B}_p}$ and contain $\theta(p,l,q)$ as its base $(n\geq p+q+l+3)$. If $plq\neq 0$, then
\begin{eqnarray*}
  i_+(G_w) &\geq& \left\{
                    \begin{array}{ll}
                     1+\frac{p+q+l}{2} , & \hbox{is $p+q+l$ is even;} \\[3pt]
                    \frac{p+q+l}{2}, & \hbox{if $p,q,l$ are odd;} \\[3pt]
                    1+\frac{p+q+l+1}{2}, & \hbox{otherwise.}
                    \end{array}
                  \right.
\end{eqnarray*}
This bound is sharp.
\end{thm}
\begin{proof}
Let $u, v$ be two vertices in $\theta(p,l,q)$ (see Fig. 1) and $\mathscr{S}_3(n)\, (n\geq p+q+l+3)$ be the set of all $n$-vertex weighted bicyclic graphs with $n-p-q-l-2$ pendant vertices attached to a vertex, different from $u$ and $v$ of $\theta(p,l,q)$. Let $G^{**}$ be the bicyclic graph with $n-p-q-l-2 \, (n\geq p+q+l+3)$ pendant vertices attached to $v$ in $\theta(p,l,q)$ and let $G^{**}_w$ denote the weighted graph with $G^{**}$ as its underlying graph, where $G^{**}$ is depicted in Fig.~2. We will verify that $i_+(G^{**}_w)\leq i_+(G_w)$ for any $G_w\in \mathscr{S}_3(n)$.

For any $G_w\in \mathscr{S}_3(n)$, without loss of generality, assume that $n-p-q-l-2$ pendant vertices are attached to a vertex of $P_{p+2}-u-v$ in $G_w$. By Lemma 2.6, we have
\begin{eqnarray*}
  i_+(G_w) &=& \left\{
               \begin{array}{ll}
                 1+\frac{p}{2}+i_+(P^{l+q+1}_w), & \hbox{if $p$ is even;} \\[3pt]
                1+\frac{p-1}{2}+i_+(C^{l+q+2}_w),\ or \
                1+\frac{p+1}{2}+i_+(P^q_w)+i_+(P^l_w), & \hbox{if $p$ is odd.}
               \end{array}
             \right.\\[3pt]
  i_+(G^{**}_w) &=& \left\{
                    \begin{array}{ll}
                      1+\frac{p}{2}+i_+(P^{l+q+1}_w), & \hbox{if $p$ is even;} \\[3pt]
                     1+\frac{p+1}{2}+i_+(P^q_w)+i_+(P^l_w), & \hbox{if $p$ is odd.}
                    \end{array}
                  \right.
\end{eqnarray*}
Note that $i_+(C^{q+l+2}_w)\geq i_+(P^q_w)+i_+(P^l_w)+1$ from Lemma 2.3, hence we have $i_+(G^{**}_w)\leq i_+(G_w)$.

By Lemmas 2.9, 2.10 and 2.11, $G^{**}_w$ attains the minimal positive index among all $n$-vertex weighted bicyclic graphs with pendant vertices containing $\theta(p,l,q)$ as its base, $n\geq p+q+l+3$.
\end{proof}

Similarly, we can have the following theorem:
\begin{thm}
Let $G_w\in {\mathscr{B}_p}$ and contain $\theta(p,l,q)$ as its base $(n\geq p+q+l+3)$. If $plq\neq 0$, then
\begin{eqnarray*}
  i_-(G_w) &\geq& \left\{
                    \begin{array}{ll}
                     1+\frac{p+q+l}{2} , & \hbox{is $p+q+l$ is even;} \\[3pt]
                    \frac{p+q+l}{2}, & \hbox{if $p,q,l$ are odd;} \\[3pt]
                    1+\frac{p+q+l+1}{2}, & \hbox{otherwise.}
                    \end{array}
                  \right.
\end{eqnarray*}
This bound is sharp.
\end{thm}
Next we consider the special case that one of $p,l,q$ is zero, Without loss of generality, we may assume $l=0$. By a similar discussion as in the proof of Theorem 3.3, we can get the following result.
\begin{thm}
Let $G_w\in {\mathscr{B}_p}$ and contain $\theta(p,0,q)$ as its base $(n\geq p+q+l+3)$. Then
\begin{eqnarray*}
 && i_+(G_w)=i_-(G_w) \geq \left\{
                 \begin{array}{ll}
                   1+\frac{p+q}{2}, & \hbox{if $p+q$ is even;} \\[3pt]
                 1+\frac{p+q+1}{2}, & \hbox{otherwise.}
                 \end{array}
               \right.
\end{eqnarray*}
This bound is sharp.
\end{thm}
By Theorems 3.5 and 3.6 we have
\begin{thm}
Let $G_w$ be a weighted $\theta$-graph of order $n$ with pendant vertices. Then $i_+(G_w)\geq 2$, $i_-(G_w)\geq 2$ and $i_0(G_w)\leq n-4$.
\end{thm}

\section{\normalsize Characterization of weighted bicyclic graphs with small positive (negative) indices}\setcounter{equation}{0}
In this section we characterize the extremal weighted bicyclic graphs with positive (resp. negative) indices $1,2.$
\begin{thm}
Let $G_w\in {\mathscr{B}}$. Then $i_+(G_w)=1$ if and only if $G_w$ is one of the following graphs: the weighted graph $\theta(1,1,1)_w$ with weighted condition $c_1a_2=a_1c_2$ and $a_2b_1=a_1b_2;$ the weighted graph $\theta(1,0,1)_w$ with weighted condition $a_2c_1=a_1c_2$.
\end{thm}
\begin{proof}
By Theorems 3.3 and 3.7, it suffices to consider the case that the weighted bicyclic graphs of order $n$ without pendant vertices. If $G_w$ is a $\infty$-graph, it contains $P^2_w\bigcup P^2_w$ as an induced subgraph, hence $i_+(G_w)\geq i_+(P^2_w\bigcup P^2_w)=2$. Then we just need to consider the case that $G_w$ is a $\theta$-graph. Without loss of generality, we assume that $l\leq p\leq q$.

If $l=0$, then we have $p+q+1\leq 3$, otherwise it contains $P^4_w$ as an induced subgraph and by Lemma 2.7, $i_+(P^4_w)=2$. Noted that $p+q\geq 2$, then the underlying graph of $G_w$ must be $\theta(1,0,1)$. Applying ECMOs to $A(G_w)$ yields $i_+(G_w)=1$ if and only if the weight of $G_w$ satisfies
$a_2c_1=a_1c_2$.

If $l>0$, then we have $p+q+2\leq 4$, otherwise it contains $C^k_w$ as an induced subgraph and $i_+(C^k_w)\geq3$, where $k\geq 5$. Noted that $p+q\geq 2$, then the underlying graph of $G_w$ must be $\theta(1,1,1)$. Applying ECMOs to $A(G_w)$ yields $i_+(G_w)=1$ if and only if the weight of $G_w$ satisfies
$a_2c_1=a_1c_2$ and $a_2b_1=a_1b_2$.
\end{proof}

Similarly, we have the following theorem:
\begin{thm}
Let $G_w\in {\mathscr{B}}$. Then $i_-(G_w)=1$ if and only if $G_w$ is the weighted graph $\theta(1,1,1)_w$ with weighted condition $c_1a_2=a_1c_2$ and $a_2b_1=a_1b_2$.
\end{thm}

\renewcommand{\arraystretch}{1.6}
\begin{table}[h!]
  \caption{The weighted condition for each $G_w\in \mathscr{B} \backslash \mathscr{B}_p$ satisfying $i_+(G_w)=2.$}\label{dd}
\begin{center}
\begin{tabular}{|c|c||c|c|}
  \hline
   weighted graph $G_w$ & weighted conditions of $G_w$ & weighted graph $G$ & weighted conditions of $G_w$ \\\hline
  $\infty(3,1,3)_w$ &   & $\theta(1,1,1)_w$ & $a_2b_1\neq a_1b_2,$ or $a_2c_1\neq a_1c_2$ \\\hline
   $\infty(3,2,3)_w$ & $4a_1a_3b_1b_3-a_2b_2c^2_1\geq 0$ & $\theta(1,0,1)_w$ &  $a_2c_1\neq a_1c_2$\\\hline
   $\infty(3,1,4)_w$& $b_1b_3=b_2b_4$ & $\theta(1,0,2)_w$ & $a_1b_2\geq c_1c_3$ \\\hline
   $\infty(4,1,4)_w$& $a_1a_3=a_2a_4, b_1b_3=b_2b_4$ & $\theta(2,0,2)_w$ & $a_2b_1c_2=a_1a_3c_2+a_2c_1c_3$ \\
  \hline
\end{tabular}
\end{center}
\end{table}
\begin{thm}
Let $G_w\in {\mathscr{B}}\backslash \mathscr{B}_p$, then $i_+(G_w)=2$ if and only if  $G_w\cong \infty(3,1,3)_w, \infty(3,2,3)_w, \infty(3,1,4)_w, \linebreak \infty(4,1,4)_w, \theta(1,1,1)_w, \theta(1,0,1)_w, \theta(1,0,2)_w,$ or $\theta(2,0,2)_w$ and the corresponding weighted conditions are as shown in Table $1$, where the empty cell means there is no correlation between the inertia index of $G_w$ and its weight set.
\end{thm}
\begin{proof}
We distinguish the following two possible cases to prove our results.
\vspace{2mm}

{\bf Case 1} $G_w$ is a weighted $\infty$-graph.
\vspace{2mm}

Note that if $G_w$ contains $P^6_w$ as an induced subgraph, then $i_+(G_w)\ge 3$. Hence, it suffices to consider that $p+l+q-4\leq 5$, i.e., $p+l+q\leq 9$. Note that $p+l+q\geq 7$, hence $7\le p+l+q\leq 9.$

If $p+l+q=7$, then $G_w$ must be $\infty(3,1,3)_w$. Applying the ECMOs to $A(G_w)$, we have $i_+(G_w)=2$ and the positive index of $G_w$ is independent of its weights.

If $p+l+q=8$, then $G_w\cong \infty(3,2,3)_w$ or, $\infty(3,1,4)_w$. Applying the ECMOs to $A(G_w)$, if $G_w\cong \infty(3,2,3)_w$, then we have $i_+(G_w)=2$ if and only if the weight of $G_w$ satisfies $4a_1a_3b_1b_3-a_2b_2c^2_1\geq 0$; if $G_w\cong \infty(3,1,4)_w$, then we have $i_+(G_w)=2$ if and only if the weight of $G_w$ satisfies $b_1b_3=b_2b_4$.

If $p+l+q=9$, then $G_w\cong \infty(3,3,3)_w, \infty(3,2,4)_w, _w\infty(3,1,5)$ or, $\infty(4,1,4)_w$. Applying the ECMOs to $A(G_w)$, if $G_w\cong\infty(4,1,4)_w$, then we have $i_+(G_w)=2$ if and only if the weight of $G_w$ satisfies $a_1a_3=a_2a_4$ and $b_1b_3=b_2b_4$; if $G_w\cong \infty(3,3,3)_w, \infty(3,2,4)_w$ or $\infty(3,1,5)_w$, then $G_w$ contains $H_w$ as its induced subgraph, where the underlying graph of $H_w$ is depicted in Fig. 3. By Lemma 2.6,  $i_+(G_w)\geqslant i_+(H_w)\geq 3.$
\begin{figure}[h!]
\begin{center}
  \psfrag{9}{$H$}
   \includegraphics[width=30mm]{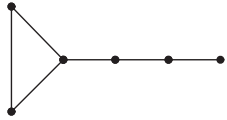}\\
  \caption{The underlying graph of $H_w$.}
\end{center}
\end{figure}

{\bf Case 2} $G_w$ is a weighted $\theta$-graph. In this case, we assume, without loss of generality, that $l\leq p\leq q$. 
By Lemmas 2.5 and 2.7, we have $i_+(P^6_w)=3$ and $i_+(C^k_w)\geq3,\, k\geq 5$. Hence, it suffices to consider that $G_w$ does not contain $P^6_w$ or $C^k_w$ as an induced subgraph, $k\geq5$.

First consider $l>0$. In this subcase, we have $p+q+2\leq 4$, otherwise $G_w$ contains $C^k_w$ as an induced subgraph with $k\geq5$. Hence, $i_+(G_w)\geq i_+(C_w^k)\ge 3.$ It is routine to check that $p+q\geq 2$, hence $ p+q =2$, which implies the underlying graph of $G_w$ must be $\theta(1,1,1)$. Applying the ECMOs to $A(G_w)$ yields $i_+(G_w)=2$ if and only if the weight of $G_w$ satisfies $a_2b_1\neq a_1b_2$ or, $a_2c_1\neq a_1c_2$.

Now consider $l=0$. In this subcase, we have $p+q+1\leq 5$; otherwise $G_w$ contains $P^6_w$ as an induced subgraph. Note that $p+q\geq 2$, hence $2\leq p+q \leq 4$.

If $p+q=2$, then $G_w\cong \theta(1,0,1)_w$. Applying ECMOs to $A(G_w)$ yields $i_+(G_w)=2$ if and only if the weight of $G_w$ satisfies $a_2c_1\neq a_1c_2$.
If $p+q=3$, then $G_w\cong \theta(1,0,2)_w$. Applying ECMOs to $A(G_w)$ yields $i_+(G_w)= 2$ if and only if the weight of $G_w$ satisfies $a_1b_2\geq c_1c_3$.
If $p+q=4$, $G_w\cong \theta(1,0,3)_w$ or, $\theta(2,0,2)_w$. If $G_w\cong \theta(2,0,2)_w$, then applying ECMOs to $A(G_w)$ yields $i_+(G_w)=2$ if and only if the weight of $G_w$ satisfies $a_2b_1c_2-a_1a_3c_2-a_2c_1c_3=0$. If $G_w\cong \theta(1,0,3)_w$, then applying ECMOs to $A(G_w)$ yields $i_+(G_w)=3$ and the positive index of $G_w$ is independent of the weights.
\end{proof}
\begin{table}[h!]
  \centering
  \caption{The weighted condition for each $G_w\in \mathscr{B}_p$ but no pendant twins and satisfying $i_+(G_w)=2.$}\label{dd}\vspace{2mm}
\begin{tabular}{|c|c||c|c|}
  \hline
  weighted graph $G_w$ & weighted conditions of $G_w$ & weighted graph $G_w$ & weighted conditions of $G_w$ \\\hline
  $G^1_w, G^3_w, G^6_w, G^7_w, G^8_w$ &   & $G^9_w, G^{10}_w$ & $a_1c_2=a_2c_1 $ \\\hline
  $G^2_w$ & $a_1b_2=a_2b_1$ & $G^{11}_w$ &  $ a_1a_3= a_2b_1$\\\hline
   $G^4_w, G^5_w$& $a_1b_2=a_2b_1$, $a_1c_2=a_2c_1$ &   &  \\\hline
\end{tabular}
\end{table}

In what follows, we shall characterize all weighted bicyclic graphs with pendants having two positive eigenvalues.
\begin{thm}
Let $G_w\in \mathscr{B}_p$ but no pendant twins. Then $i_+(G_w)=2$ if and only if $G_w\cong G^1_w, G^2_w, \ldots, G^{10}_w$ or, $G^{11}_w$ and the corresponding weighted conditions are as shown in Table $2$, where the underlying graphs of $G^1_w, G^2_w, \ldots, G^{10}_w, G^{11}_w$ are depicted in Fig. 4 and  the empty cell in Table 2 means there is no correlation between the inertia index of $G_w$ and its weight set.
\end{thm}
\begin{figure}[h!]
\begin{center}
  \psfrag{1}{$G^{1}$}\psfrag{6}{$G^{6}$}
  \psfrag{2}{$G^{2}$}\psfrag{7}{$G^{7}$}
  \psfrag{3}{$G^{3}$}\psfrag{8}{$G^{8}$}
  \psfrag{4}{$G^{4}$}\psfrag{9}{$G^{9}$}
  \psfrag{5}{$G^{5}$}\psfrag{a}{$G^{10}$}
  \psfrag{b}{$G^{11}$}
  \includegraphics[width=110mm]{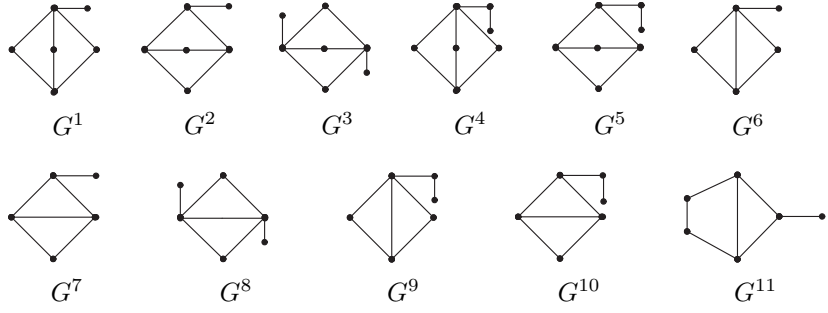}\\
  \caption{Graphs $G^1, G^2, \ldots, G^{11}$.}
\end{center}
\end{figure}
\begin{figure}[h!]
\begin{center}
  \psfrag{1}{$G^{12}$}\psfrag{6}{$G^{17}$}
  \psfrag{2}{$G^{13}$}\psfrag{7}{$G^{18}$}
  \psfrag{3}{$G^{14}$}\psfrag{8}{$G^{19}$}
  \psfrag{4}{$G^{15}$}\psfrag{9}{$G^{20}$}
  \psfrag{5}{$G^{16}$}\psfrag{a}{$G^{21}$}
  \psfrag{b}{$G^{22}$}\psfrag{c}{$G^{23}$}
  \psfrag{d}{$G^{24}$}\psfrag{e}{$G^{25}$}
  \psfrag{f}{$G^{26}$}\psfrag{g}{$G^{27}$}
  \psfrag{h}{$G^{28}$}\psfrag{i}{$G^{29}$}
  \psfrag{j}{$G^{30}$}\psfrag{k}{$G^{31}$}
  \psfrag{l}{$G^{32}$}\psfrag{m}{$G^{33}$}
  \psfrag{n}{$G^{34}$}\psfrag{p}{$G^{35}$}\psfrag{q}{$G^{36}$}
    \includegraphics[width=130mm]{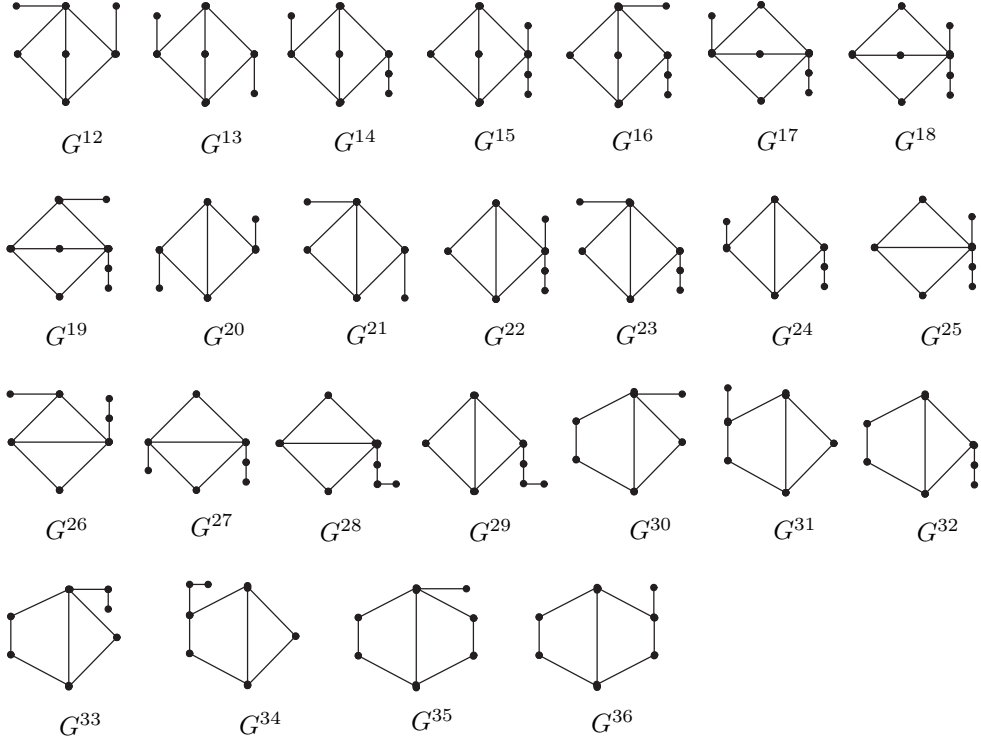}\\
  \caption{Graphs $G^{12}, G^{13}, \ldots, G^{36}.$}
\end{center}
\end{figure}
\begin{proof}By Lemmas 2.5, 2.6 and applying the ECMOs, it is routine to check that $i_+(G^i_w)=2,\,i=1,\ldots, 11$, and the weight condition for $G^i_w$ is listed in Table 2. Furthermore, $i_+(G^j_w)\geq3$ holds for any weighted condition, $j=12,\ldots, 36$. Here the underlying graphs of $G^1_w, \ldots, G^{11}_w$ are depicted in Fig. 4, while those of $G^{12}_w, \ldots, G^{36}_w$ are depicted in Fig. 5.

Let $\mathcal{H}=G-\chi(G)$ and denote by $v(\mathcal{H})$ the number of vertices of $\mathcal{H}$ in what follows.
Note that if $i_+(G_w)=2$, by Theorem 3.3, $G_w$ must be a weighted $\theta$-graph and by Lemma 2.4, we have $i_+(\chi(G_w))\leq 2$. Hence, in view of Theorems 4.1 and 4.3, we have $\chi(G)\in \{\theta(1,1,1),\theta(1,0,1),\theta(1,0,2), \theta(2,0,2)\}$.

First, we characterize all graphs $G_w$ with $\theta(1,1,1)$ as its base satisfying $i_+(G_w) = 2$ according to the following two possible cases.\vspace{2mm}

{\bf Case 1.}\ $\mathcal{H}$ is a collection of isolated vertices.\vspace{2mm}

If $v(\mathcal{H})=1$, $G_w$ must be $G^{1}_w$ or $G^{2}_w$. By Lemmas 2.5 and 2.6, we have $i_+(G^{2}_w)=2$ if and only if  the weight of $G_w$ satisfies  $a_1b_2=a_2b_1$. It is routine to check that $i_+(G^{1}_w)=2$.

If $v(\mathcal{H})=2$, $G_w$ must be $G^{3}_w, G^{12}_w$ or, $G^{13}_w$, but $i_+(G^{12}_w)=i_+(G^{13}_w)=3$.

If $v(\mathcal{H})\geq 3$, then by Lemma 2.4 $i_+(G_w)\geq 3$ since $G_w$ contains $G^{12}_w$ or, $G^{13}_w$ as an induced subgraph.\vspace{2mm}

{\bf Case 2.}\ $\mathcal{H}$ has a $P_2$ as an induced subgraph.\vspace{2mm}

If $\mathcal{H}=P_2$, $G_w$ must be $G^{4}_w$ or $G^{5}_w$. By Lemma 2.6, $i_+(G^{4}_w)=i_+(G^{5}_w)=1+i_+(G'_w)$, where $G'_w$ is $\theta (1,1,1)_w$. By applying ECMOs on $A(G'_w)$, we have $i_+(G^{4}_w)=i_+(G^{5}_w)= 2$ if and only if the weight of $G_w$ satisfies the condition that $a_2b_1=a_1b_2$ and $a_2c_1=a_1c_2$.

If $\mathcal{H}$ contains the union of $P_2$ and an isolated vertex as an induced subgraph, then by Lemma 2.4, $i_+(G_w)\geq 3$ since it contain one of $G^i_w$'s ($i=14,\ldots ,19$) as an induced subgraph.

If $\mathcal{H}$ contains a $P_3$ as an induced subgraph, then by Lemma 2.4, $i_+(G_w)=1+i_+(G'_w)\geq 3$, where $G'_w$ is $G^{1}_w$ or $G^{2}_w$ and $i_+(G^1_w)=2, i_+(G^2_w)\geq 2$.

Next we characterize all graphs $G_w$ with $\theta(1,0,1)$ as its base satisfying $i_+(G_w) = 2$ according to the following four possible cases.\vspace{2mm}

{\bf Case 1.} $\mathcal{H}$ is a collection of isolated vertices.\vspace{2mm}

If $v(\mathcal{H})=1$, $G_w$ must be $G^{6}_w$ or $G^{7}_w$.

If $v(\mathcal{H})=2$, $G_w$ must be $G^{8}_w, G^{20}_w$ or $G^{21}_w$, but $i_+(G^{20}_w)=i_+(G^{21}_w)=3$.

If $v(\mathcal{H})\geq 3$, then by Lemma 2.4, $i_+(G_w)\geq 3$ since $G_w$ contains $G^{20}_w$ or $G^{21}_w$ as an induced subgraph.\vspace{2mm}

{\bf Case 2.} $\mathcal{H}$ is $P_2$. In this subcase, the underlying graph of $G_w$ must be $G^{9}_w$ or $G^{10}_w$, by calculation we have $i_+(G^{9}_w)=i_+(G^{10}_w)=2$ if and only if the weight of $G^{9}_w$ and $G^{10}_w$ satisfies the condition that $a_2c_1=a_1c_2$.\vspace{2mm}

{\bf Case 3.} $\mathcal{H}$ contains the union of $P_2$ and an isolated vertex as an induced subgraph. By Lemma 2.4, $i_+(G_w)\geq 3$ since it contains one of $G^i_w$'s ($i=22,\ldots ,27$) as an induced subgraph.\vspace{2mm}

{\bf Case 4.} $\mathcal{H}$ contains a $P_3$ as an induced subgraph. By Lemma 2.4, $i_+(G_w)\geq 3$ since $G_w$ contains $G^{28}_w$ or $G^{29}_w$ as an induced subgraph.

Now we characterize all graphs $G_w$ with $\theta(1,0,2)$ as its base satisfying $i_+(G_w) = 2$ according to the following two possible cases.\vspace{2mm}

{\bf Case 1.} $\mathcal{H}$ is a collection of isolated vertices.\vspace{2mm}

If $v(\mathcal{H})=1$, $G_w$ must be $G^{11}_w, G^{30}_w$ or $G^{31}_w.$ Note that $i_+(G^{30}_w)=i_+(G^{31}_w)=3$, and by Lemmas 2.5 and 2.6, $i_+(G^{11}_w)=2$ if and only if the weight of $G_w^{11}$ satisfies $a_1a_3=a_2b_1$.

If $v(\mathcal{H})\geq 2$, then by Lemma 2.4, $i_+(G_w)\geq 3$ since $G_w$ contains $G^{30}_w$ or $G^{31}_w$ as an induced subgraph.\vspace{2mm}

{\bf Case 2.} $\mathcal{H}$ contains a $P_2$ as a induced subgraph, $G_w$ must be $G^{32}_w, G^{33}_w$ or $G^{34}_w$, but each of them have more than 2 positive eigenvalues.

At last, we consider graphs $G_w$ with $\theta(2,0,2)$ as its base satisfying $i_+(G_w) = 2.$ In fact, in this case, $G_w$ contains $G^{35}_w$ or $G^{36}_w$ as an induced subgraph.
\end{proof}

Similarly, we can have the following theorems.
\begin{thm}
Let $G_w\in {\mathscr{B}}\backslash \mathscr{B}_p$, then $i_-(G_w)=2$ if and only if $G_w$ is one of the following graphs: the weighted graph $\infty(4,1,4)_w$ with weighted condition $a_1a_3=a_2a_4$ and $b_1b_3=b_2b_4$; the weighted graph $\theta(1,1,1)_w$ with weighted condition $a_2b_1\neq a_1b_2$ or $a_2c_1\neq a_1c_2$; the weighted graph $\theta(1,0,1)_w$; the weighted graph $\theta(1,0,2)_w$ with weighted condition $a_1b_2\leq c_1c_3$; the weighted graph $\theta(2,0,2)_w$ with weighted condition $a_2b_1c_3-a_1a_3c_2-a_2c_1c_3=0$; the weighted graph $\theta(1,1,2)_w$ with weighted condition $a_1b_2=a_2b_1$.
\end{thm}
\begin{thm}
Let $G_w\in \mathscr{B}_p$ but no pendant twins, then $i_-(G_w)=2$ if and only if $G_w$ is one of the following graphs: the weighted graph $G^{1}_w, G^{3}_w,G^{6}_w, G^{8}_w$; the weighted graph $G^{4}_w, G^{5}_w$) with weighted condition $a_2b_1=a_1b_2$ and $a_2c_1=a_1c_2$; the weighted graph $G^{2}_w$ with weighted condition $a_1b_2=a_2b_1$; the weighted graph $G^{11}_w$ with weighted condition $a_1a_3=a_2b_1.$
\end{thm}

\section{\normalsize Weighted bicyclic graphs with rank $2,3,4$}\setcounter{equation}{0}
The rank of a weighted bicyclic graph $G_w$ is the rank of its adjacency matrix $A(G_w)$, denoted by $r(G_w)$. Then it is easy to see that $r(G_w)=i_+(G_w)+i_-(G_w)$. In this section, we'll characterize the weighted bicyclic graphs with rank $2,3,4,$ respectively.
\begin{thm}\label{thm5.1}
Let $G_w\in {\mathscr{B}}$, then $r(G_w)=2$ if and only if $G_w\cong \theta(1,1,1)_w$ with weighted condition $a_1c_2=a_2c_1$ and $a_1b_2=a_2b_1$.
\end{thm}
\begin{proof}
Let $G_w$ be a weighted bicyclic graph, $i_+(G_w)\geq 1$ and $i_-(G_w)\geq 1$ since $G$ contains $P_2$ as an induced subgraph. Then $r(G_w)=2$ if and only if $i_+(G_w)=i_-(G_w)=1$. By Theorems 4.1 and 4.2 we know $G_w$ must be $\theta(1,1,1)_w$ with weighted condition $a_1c_2=a_2c_1$ and $a_2b_1=a_1b_2$.
\end{proof}

\begin{thm}
Let $G_w\in {\mathscr{B}}$, then $r(G_w)=3$ if and only if $G_w\cong \theta(1,0,1)_w$ with weighted condition $a_2c_1=a_1c_2$.
\end{thm}
\begin{proof}
Let $G_w$ be a weighted bicyclic graph, since $i_+(G_w)\geq 1$ and $i_-(G_w)\geq 1$, then $r(G_w)=3$ if and only if $i_+(G_w)=1,i_-(G_w)=2$ or $i_+(G_w)=2,i_-(G_w)=1$. Note that either $i_+(G_w)$ or $i_-(G_w)$ equals 1, hence by Theorems 4.1 and 4.2 we know $G_w$ must be $\theta(1,0,1)_w$ with weighted condition $a_2c_1=a_1c_2$.
\end{proof}

\begin{thm}
Let $G_w\in {\mathscr{B}}\backslash \mathscr{B}_p$, then $r(G_w)=4$ if and only if $G_w$ is one of the following graphs: the weighted graph $\infty (4,1,4)_w$ satisfying $a_1a_3=a_2a_4$ and $b_1b_3=b_2b_4$; the weighted graph $\theta (1,1,1)_w$ satisfying $a_2b_1\neq a_1b_2$ or $a_2c_1\neq a_1c_2$; the weighted graph $\theta (1,0,1)_w$ satisfying $a_2c_1\neq a_1c_2$; the weighted graph $\theta (1,0,2)_w$ satisfying $a_1b_2=c_1c_3$; the weighted graph $\theta (2,0,2)_w$ satisfying $a_2b_1c_3-a_1a_3c_2-a_2c_1c_3=0$.
\end{thm}
\begin{proof}
If $G_w$ be a weighted bicyclic graph, it is easy to know that $i_+(G_w)\geq 1$ and $i_-(G_w)\geq 1$. Then $r(G_w)=4$ if and only if $(i_+(G_w),i_-(G_w))=(1,3)$ or $(i_+(G_w),i_-(G_w))=(3,1)$ or $(i_+(G_w),i_-(G_w))=(2,2)$. If one of $i_+(G_w)$ and $i_-(G_w)$ equals 1, then $G_w$ must be $\theta(1,1,1)_w$ or $\theta(1,0,1)_w$, by Theorems 4.1 and  4.2 we know $r(G_w)<4$.

Hence, it suffices to consider that $(i_+(G_w),i_-(G_w))=(2,2)$. By Theorems 4.3 and 4.5, $(i_+(G_w),i_-(G_w))=(2,2)$ if and only if $G_w$ is one of the graphs described in Theorem 5.3.
\end{proof}

Similarly, we can have the following theorem:
\begin{thm}
Let $G_w\in \mathscr{B}_p$ but no pedant twins $n(n\geq 4)$, if $r(G_w)=4$ if and only if $G_w$ is one of the following graphs: the weighted graphs $G^{1}_w, G^{3}_w, G^{6}_w, G^{8}_w$; the weighted graph $G^{2}_w$ satisfying the weighted condition $a_1b_2=a_2b_1$; the weighted graph $G^{4}_w, G^{5}_w$ satisfying the weighted condition $a_2b_1=a_1b_2$ and $a_2c_1=a_1c_2$; the weighted graph $G^{11}_w$ satisfying the weighted condition $a_1a_3=a_2b_1$.
\end{thm}

\end{document}